# FEEDBACK STABILIZATION METHODS FOR THE SOLUTION OF NONLINEAR PROGRAMMING PROBLEMS

**Iasson Karafyllis**

**Department of Environmental Engineering,**
Technical University of Crete, 73100, Chania, Greece
email: ikarafyl@enveng.tuc.gr

**Abstract**

In this work we show that given a nonlinear programming problem, it is possible to construct a family of dynamical systems defined on the feasible set of the given problem, so that: (a) the equilibrium points are the unknown critical points of the problem, (b) each dynamical system admits the objective function of the problem as a Lyapunov function, and (c) explicit formulae are available without involving the unknown critical points of the problem. The construction of the family of dynamical systems is based on the Control Lyapunov Function methodology, which is used in mathematical control theory for the construction of stabilizing feedback. The knowledge of a dynamical system with the previously mentioned properties allows the construction of algorithms which guarantee global convergence to the set of the critical points.

## 1. Introduction

Differential equations have been used in the past for the solution of Nonlinear Programming (NLP) problems. The reader may consult [1,6,7,11,21,25,26,27,28] for various results on the topic. Some methods are interior-point methods (in the sense that are defined only on the feasible set) while other methods are exterior-point methods (in the sense that are defined at least in a neighborhood of the feasible set). As remarked in [5], each system of differential equations that solves a NLP problem when combined with a numerical scheme for solving Ordinary Differential Equations (ODEs) provides a numerical scheme for solving the NLP problem.

In this work, we are interested in the application of feedback stabilization methods for solving NLP problems. The feedback stabilization methods can be applied in two ways:
> ➢ first, for the construction of the dynamical system, which solves the NLP problem and
> ➢ for the selection of the step size of the Runge-Kutta scheme that is used for the solution of the resulting system of ODEs (see [12,16]).

More specifically, consider the Nonlinear Programming problem:

$$\min\{\theta(x) : x \in S\} \quad (1.1)$$

where $x \in \Re^n$ and the closed set $S \subseteq \Re^n$ is defined by

$$S := \left\{ x \in \Re^n : h_1(x) = ... = h_m(x) = 0, \max_{j=1,...,k}(g_j(x)) \leq 0 \right\} \quad (1.2)$$

where $m < n$ and all functions $\theta : \Re^n \to \Re$, $h_i : \Re^n \to \Re$ ($i=1,...,m$), $g_j : \Re^n \to \Re$ ($j=1,...,k$) are twice continuously differentiable. Inspired by the methods employed in the book [15], we would like to construct a well-defined dynamical system $\dot{x} = F(x)$ on the closed set $S \subseteq \Re^n$, where $F : S \to \Re^n$ is a continuous vector field with the following properties:

Property 1: For every $x \in S$, $F(x)$ belongs to the contingent cone to $S$ at $x$. This property is required because local existence of solutions of the dynamical system $\dot{x} = F(x)$ can be guaranteed by Nagumo's theorem (given on page 27 of the book [3]).

Property 2: $F: S \to \Re^n$ is a locally Lipschitz vector field. This property is required for uniqueness of the solutions of the dynamical system $\dot{x} = F(x)$. Moreover, this property is required because we would like to be able to apply $1^{st}$ order Runge-Kutta schemes for the simulation of the solutions of the dynamical system $\dot{x} = F(x)$. Higher regularity is also desirable because high order Runge-Kutta schemes can be used for the simulation of the solutions of the dynamical system $\dot{x} = F(x)$.

Property 3: The equilibrium points of the dynamical system $\dot{x} = F(x)$ are exactly the points $x \in S$ for which there exist $\lambda_i \in \Re$ ($i = 1,...,m$) and $\mu_j \geq 0$ ($j = 1,...,k$) such that the necessary Karush-Kuhn-Tucker conditions hold:

$$\nabla \theta(x) + \sum_{i=1}^{m} \lambda_i \nabla h_i(x) + \sum_{j=1}^{k} \mu_j \nabla g_j(x) = 0$$
$$\sum_{j=1}^{k} \mu_j g_j(x) = 0 \tag{1.3}$$

Property 4: The function $V(x) = \theta(x) - \theta(x^*)$ is a Lyapunov function for the dynamical system $\dot{x} = F(x)$, where $x^* \in S$ is one of the global solutions of the NLP described by (1.1), (1.2), i.e., $\theta(x^*) = \min\{\theta(x): x \in S\}$. In other words, we would like the inequality $\nabla \theta(x) F(x) < 0$ to hold for all $x \in S$ for which there are no $\lambda_i \in \Re$ ($i = 1,...,m$) and $\mu_j \geq 0$ ($j = 1,...,k$) such that conditions hold (1.3). This property is important because it guarantees useful stability properties. Furthermore, the fact that the Lyapunov function of the dynamical system $\dot{x} = F(x)$ has the special form $V(x) = \theta(x) - \theta(x^*)$ is important for numerical purposes (see [12,16]): the time derivative of the Lyapunov function along the solutions of the system and the difference of the values of the Lyapunov function between two points can be computed without knowledge of the solution $x^* \in S$ of the NLP problem.

Property 5: The vector field $F: S \to \Re^n$ must be explicitly known. Formulae for the vector field $F: S \to \Re^n$ must be provided: the formulae must not involve the solution $x^* \in S$ of the NLP described by (1.1), (1.2).

Property 6: The vector field $F: S \to \Re^n$ must have free parameters which can be selected in an appropriate way so that the convergence properties of the corresponding numerical schemes to the global attractor of the dynamical system are optimal. In other words, we want to construct a family of vector fields $F: S \to \Re^n$ with all the above properties.

It must be noted that the properties 1-6 are rarely satisfied by other differential equation methods for solving NLPs. For example, in [1] and [6], the constructed Lyapunov function is $V(x) = \frac{1}{2}|x - x^*|^2$ and this does not meet our requirements. Moreover, for [6] the point $x = x^*$ is not an equilibrium point for the constructed time-varying system $\dot{x} = F(t, x)$. Antipin in [1] constructs an autonomous system $\dot{x} = F(x)$ for which $x = x^*$ is an equilibrium point and $F \in C^0(\Re^n; \Re^n)$ (being locally Lipschitz) does not depend on the location of the unknown point $x^* \in \Re^n$. However, the computation of $F(x)$ is involved (it requires the solution of a (NLP) since it involves a projection on the feasible set). The (NLP) without equality constraints under additional convexity hypotheses has been studied in [26]. However, again the constructed Lyapunov function is of the form $V(x) = P(x) + \frac{1}{2}|x - x^*|^2$ and this does not meet our requirements. On the other hand, the papers [25,27] propose systems of differential equations that satisfy properties 1-6 for systems without inequality constraints. Local results are provided in the paper [28] and differential equations based on barrier methods were considered in [7].

It is clear that the knowledge of the Lyapunov function $V(x) = \theta(x) - \theta(x^*)$ can allow us to construct the vector field $F : S \to \Re^n$ by the Control Lyapunov Function methodology of feedback design (see [2,9,17,23]) for the control system $\dot{x} = u \in \Re^n$. However, there are certain obstructions for the direct application of the classical Control Lyapunov Function methodology: (i) the system is not defined on $\Re^n$ but on the closed set $S \subseteq \Re^n$, (ii) for every $x \in S$, $F(x)$ must belong to the contingent cone to $S$ at $x$, and (iii) the position of the equilibrium points, i.e., the set of points $x \in S$ for which there exist $\lambda_i \in \Re$ ($i = 1,...,m$) and $\mu_j \geq 0$ ($j = 1,...,k$) such that conditions (1.3) hold is unknown (this is what we are looking for).

The contribution of the paper is twofold:

- The main result of the present work (Theorem 2.1) shows that all the previously mentioned obstructions can be overcome under appropriate assumptions.

- Based on the ideas described in [12,16], in Section 3 of the present work, we present an algorithm for the solution of the NLP described by (1.1) and (1.2) which is based on the application of the explicit Euler scheme for the numerical solution of the resulting system of ODEs with appropriate step selection (Theorem 3.1). The algorithm will converge for every initial condition (global convergence). A modified and simpler version of the algorithm can work under slightly more demanding assumptions (Remark 3.5).

It should be noticed that the convergence rates of the proposed algorithms depend on the selection of certain matrices which are the free parameters described in Property 5 above. However, since the proposed algorithms are global, it can be used in combination with any other local algorithm that guarantees fast convergence based on the following intuitive idea: "apply the newly proposed algorithms when you are away from a solution and apply a fast local algorithm when you are close to a solution".

It should be emphasized that no claim is made about the effectiveness of the proposed algorithms. The topic of the numerical solution of NLPs is a mature topic and it is clear that other algorithms have much better characteristics than the algorithms proposed in this paper. However, the theory used for the construction of the algorithm is different from other existing algorithms. The algorithms contained in this work are derived by using concepts of dynamical systems theory and mathematical control theory.

The structure of the paper is as follows: Section 2 contains the statement and proof of Theorem 2.3, which provides the solution to the problem of the construction of a vector field with properties 1-6. Section 3 provides numerical algorithms for the exploitation of the constructed vector field. Section 4 of the paper provides some examples, which show the performance of the algorithms. Finally, Section 5 of the paper contains the concluding remarks. The appendix provides the proofs of certain auxiliary results.

**Notations** Throughout this paper we adopt the following notations:

∗ Let $A \subseteq \Re^n$ be a set. By $C^0(A;\Omega)$, we denote the class of continuous functions on $A$, which take values in $\Omega$. By $C^k(A;\Omega)$, where $k \geq 1$ is an integer, we denote the class of differentiable functions on $A$ with continuous derivatives up to order $k$, which take values in $\Omega$. By $C^\infty(A;\Omega)$, we denote the class of differentiable functions on $A$ having continuous derivatives of all orders (smooth functions), which take values in $\Omega$, i.e., $C^\infty(A;\Omega) = \bigcap_{k \geq 1} C^k(A;\Omega)$.

∗ For a vector $x \in \Re^n$ we denote by $|x|$ its usual Euclidean norm and by $x'$ its transpose. For a real matrix $A \in \Re^{n \times m}$ we denote by $|A|$ its induced norm, i.e., $|A| := \max\{|Ax| : x \in \Re^m, |x| = 1\}$ and by $A' \in \Re^{m \times n}$ its transpose. $I_n \in \Re^{n \times n}$ denotes the identity matrix. For every $x = (x_1,...,x_n)' \in \Re^n$ we define $x^+ = (\max(0, x_1),..., \max(0, x_n))' \in \Re^n$. Notice that the following property holds for every positive definite and diagonal matrix $R \in \Re^{n \times n}$: $x'Rx^+ = 0 \Leftrightarrow x^+ = 0$.

∗ $\Re_+^n := (\Re_+)^n = \{(x_1,...,x_n)' \in \Re^n : x_1 \geq 0,...,x_n \geq 0\}$. Let $x, y \in \Re^n$. We say that $x \leq y$ if and only if $(y - x) \in \Re_+^n$.

∗ For every scalar continuously differentiable function $V : \Re^n \to \Re$, $\nabla V(x)$ denotes the gradient of $V$ at $x \in \Re^n$, i.e., $\nabla V(x) = \left(\frac{\partial V}{\partial x_1}(x),...,\frac{\partial V}{\partial x_n}(x)\right)$ and $\nabla^2 V(x)$ denotes the Hessian matrix of $V$ at $x \in \Re^n$.

## 2. Transforming an NLP problem into a feedback stabilization problem

Consider the NLP problem described by (1.1) and (1.2) under the following assumptions:

**(H1)** *The feasible set $S \subseteq \Re^n$ defined by (1.2) is non-empty and the level sets of $\theta : \Re^n \to \Re$ are compact sets, i.e., for every $x_0 \in S$ the level set*

$$\{x \in S : \theta(x) \leq \theta(x_0)\}$$

*is compact.*

**(H2)** *For every $x \in S$ the row vectors $\nabla h_i(x)$ ($i = 1,...,m$) and $\nabla g_j(x)$ for all $j = 1,...,k$ for which $g_j(x) = 0$ (active constraints) are linearly independent.*

Assumption (H1) is a standard assumption which guarantees that the NLP problem described by (1.1) and (1.2) is well-posed and admits at least one global solution (see [4]). Assumption (H2) is an extension of the main assumption in [20]. Assumption (H2) in conjunction with the main result in [20] guarantees that for every solution $x^* \in S$ of the NLP problem described by (1.1) and (1.2), there exist $\lambda_i^* \in \Re$ ($i = 1,...,m$) and $\mu_j^* \geq 0$ ($j = 1,...,k$) such that conditions (1.3) hold.

We define:

$$h(x) = \begin{bmatrix} h_1(x) \\ \vdots \\ h_m(x) \end{bmatrix} \in \Re^m, \; A(x) = \begin{bmatrix} \nabla h_1(x) \\ \vdots \\ \nabla h_m(x) \end{bmatrix} \in \Re^{m \times n}, \; g(x) = \begin{bmatrix} g_1(x) \\ \vdots \\ g_k(x) \end{bmatrix} \in \Re^k, \; B(x) = \begin{bmatrix} \nabla g_1(x) \\ \vdots \\ \nabla g_k(x) \end{bmatrix} \in \Re^{k \times n}, \text{ for all } x \in \Re^n \quad (2.1)$$

Assumption (H2) allows us to define the symmetric matrix:

$$H(x) = I_n - A'(x)(A(x)A'(x))^{-1} A(x), \text{ for all } x \in \Re^n \text{ in a neighborhood of } S \qquad (2.2)$$

The following facts are direct consequences of definition (2.2):

**Fact 1:** $H^2(x) = H(x)$, $A(x)H(x) = 0$ and $H(x)A'(x) = 0$.

**Fact 2:** $\xi' H(x)\xi = |H(x)\xi|^2$, for all $\xi \in \Re^n$

**Fact 3:** For every $\xi \in \Re^n$ there exists $\lambda \in \Re^m$ such that $\xi = H(x)\xi + A'(x)\lambda$.

Next, we define the set of critical points for the NLP problem defined by (1.1) and (1.2).

**Definition 2.1:** *Let $\Phi \subseteq S$ be the set of all points $x \in S$ for which there exist $\lambda \in \Re^m$ and $\mu \in \Re_+^k$ such that conditions (1.3) hold. In other words, $\Phi \subseteq S$ is the set of critical points or Karush-Kuhn-Tucker points for the problem defined by (1.1) and (1.2).*

Clearly, assumptions (H1) and (H2) guarantee that the set $\Phi \subseteq S$ is non-empty.

The following lemma provides a useful consequence of assumption (H2). Its proof is provided at the Appendix.

**Lemma 2.2:** *If assumption (H2) holds then the matrix*

$$Q(x) := B(x)H(x)B'(x) - diag(g(x)) \qquad (2.3)$$

*is positive definite for all $x \in S$.*

We are now ready to state the main result of this section.

**Theorem 2.3:** *Suppose that assumptions (H1) and (H2) hold for the NLP problem described by (1.1) and (1.2). Let $Q(x) \in \Re^{k \times k}$ be the symmetric positive definite matrix defined by (2.3). Let $R_1(x) \in \Re^{n \times n}$ be an arbitrary $C^1$, symmetric and positive definite matrix, $R_2(x) \in \Re^{k \times k}$ be an arbitrary $C^1$, symmetric and positive semidefinite matrix, $a_i(x)$, $b_i(x)$, $c_i(x)$ ($i = 1,...,k$) be arbitrary $C^1$ non-negative functions with $b_i(x) + c_i(x) > 0$ for all $i = 1,...,k$, $x \in S$ and at least one of the matrices $R_2(x) \in \Re^{k \times k}$, $diag(a(x)) \in \Re^{k \times k}$ being positive definite, where $a(x) := (a_1(x),...,a_k(x))' \in \Re^k$.*

*Define the following locally Lipschitz vector field:*

$$F(x) = -[H(x) - P'(x)Q(x)P(x)]R_1(x)[H(x) - P'(x)Q(x)P(x)](\nabla \theta(x))'$$
$$- P'(x)diag(g(x))(R_2(x)diag(g(x)) - diag(a(x)))v(x) \quad (2.4)$$
$$- P'(x)R_3(x)(v(x))^+$$

*where*

$$P(x) := Q^{-1}(x)B(x)H(x) \in \Re^{k \times n}$$
$$v(x) := P(x)(\nabla \theta(x))' \in \Re^k \quad (2.5)$$
$$R_3(x) := diag\left(b_1(x) + c_1(x)(\max(0, v_1(x)))^{2p_1},...,b_k(x) + c_k(x)(\max(0, v_k(x)))^{2p_k}\right)$$

*and $p_i \geq 1$ ($i = 1,...,k$) are integers.*

*Then the following properties hold:*
*a) $A(x)F(x) = 0$ for all $x \in S$,*
*b) $\nabla \theta(x)F(x) < 0$, for all $x \in S \setminus \Phi$,*
*c) $F(x) = 0 \Leftrightarrow x \in \Phi$*
*d) $B(x)F(x) = diag(g(x))Q^{-1}(x)w(x) - R_3(x)(v(x))^+$, for all $x \in S$, where*

$$w(x) := B(x)H(x)R_1(x)[H(x) - P'(x)Q(x)P(x)](\nabla \theta(x))'$$
$$- [Q(x) + diag(g(x))](R_2(x)diag(g(x)) - diag(a(x)))v(x) - R_3(x)(v(x))^+ .$$

*Consider the dynamical system*

$$\dot{x} = F(x) \quad (2.6)$$

*on the closed set $S \subseteq \Re^n$. Then the following properties hold:*

1) *For every $x_0 \in S$ there exists a unique solution $x(t)$ of the initial value problem (2.6) with $x(0) = x_0$ which is defined for all $t \geq 0$ and satisfies $x(t) \in S$ for all $t \geq 0$.*

2) *Every point $x \in \Phi$ is an equilibrium point for (2.6). Every strict local solution $x^* \in S$ of the NLP problem described by (1.1) and (1.2) is locally asymptotically stable for system (2.6).*

*If we further denote by $\omega(x_0)$ the set of accumulation points of the set $\{x(t): t \geq 0\}$, where $x_0 \in S$, then it holds that $\omega(x_0)$ is a compact, positively invariant set for which there exists $l \leq \theta(x_0)$ such that $\omega(x_0) \subseteq \Phi \cap \{x \in S : \theta(x) = l\}$.*

**Remark 2.4:** Clearly, the matrices $R_1(x) \in \Re^{n \times n}$, $R_2(x) \in \Re^{k \times k}$ and the functions $a_i(x)$, $b_i(x)$, $c_i(x)$ ($i = 1,...,k$), can be selected in an appropriate way so that the convergence properties of the corresponding numerical schemes to the global attractor of the dynamical system are optimal. The stability properties of system (2.6) are analogous to the stability properties of gradient systems (see [24]).

**Remark 2.5:** It should be noted that all properties 1-6 mentioned in the Introduction are satisfied for the dynamical system (2.6). Indeed,
-- Property 1 is a direct consequence of (a) and (d). More specifically, since $R_3(x) := diag\left(b_1(x) + c_1(x)(\max(0, v_1(x)))^{2p_1},...,b_k(x) + c_k(x)(\max(0, v_k(x)))^{2p_k}\right)$ for certain non-negative functions $b_j(x), c_j(x)$ ($j=1,...,k$) and since $\frac{d}{dt}g(x) = B(x)\dot{x} = B(x)F(x) = diag(g(x))Q^{-1}(x)w(x) - R_3(x)(v(x))^+$ it follows that the following implication holds: "if $g_j(x) = 0$ for some $j \in \{1,...,k\}$ then $\frac{d}{dt}g_j(x) = -\left(b_j(x) + c_j(x)(\max(0, v_j(x)))^{2p_j}\right)\max(0, v_j(x)) \leq 0$". The previous implication and property (a) guarantee that for every $x \in S$, $F(x)$ belongs to the contingent cone to $S$ at $x$.

-- Property 2 is a direct consequence of definitions (2.2), (2.3), (2.4), (2.5) and the fact that all functions $\theta: \Re^n \to \Re$, $h_i: \Re^n \to \Re$ ($i=1,...,m$), $g_j: \Re^n \to \Re$ ($j=1,...,k$) are twice continuously differentiable. It should be noticed that if at least one of the functions $b_i(x)$, ($i=1,...,k$) takes positive values then the vector field $F(x)$ defined by (2.4) is simply locally Lipschitz and not $C^1$. When $b_i(x) \equiv 0$, for $i=1,...,k$ then the vector field $F(x)$ defined by (2.4) is $C^1$. Higher regularity is possible by assuming higher regularity for all functions and matrices involved in (2.2), (2.3), (2.4), (2.5), sufficiently large values for the integers $p_i \geq 1$ ($i=1,...,k$) and $b_i(x) \equiv 0$, for $i=1,...,k$.

--Property 3 is a direct consequence of (c). Property 4 is a direct consequence of (b). Indeed, notice that the function $V(x) = \theta(x) - \theta(x^*)$, where $x^* \in S$ is one of the global solutions of the NLP described by (1.1), (1.2) satisfies $\nabla V(x)F(x) = \nabla \theta(x)F(x)$.

--Finally, properties 5 and 6 are evident.

**Remark 2.6:** The inspiration for Theorem 2.3 is the transformation of the NLP problem into a feedback stabilization problem. First, we notice that the Control Lyapunov Function (see [2,9,17,23]) is selected to be $V(x) = \theta(x) - \theta(x^*)$, where $x^* \in S$ is one of the global solutions of the NLP problem described by (1.1), (1.2). The only problem is that we must define in an appropriate way the control system so that $S$ is a positively invariant set for all possible inputs. In other words, we must have:

$$\frac{d}{dt}h(x) = A(x)\dot{x} = 0 \text{ and } \frac{d}{dt}g(x) = B(x)\dot{x} = diag(g(x))v - u$$

for all possible inputs $v, u \in \Re^n$. Notice that the property $\frac{d}{dt}g(x) = B(x)\dot{x} = diag(g(x))v - u$ for arbitrary $v \in \Re^n, u \in \Re^n_+$ guarantees the implication: "if $g_j(x) = 0$ for some $j \in \{1,...,k\}$ then $\frac{d}{dt}g_j(x) = -u_j \leq 0$". The property $\frac{d}{dt}h(x) = A(x)\dot{x} = 0$ implies that $\dot{x} = H(x)w$, for arbitrary $w \in \Re^n$. Combining, we get $B(x)H(x)w = diag(g(x))v - u$. By redefining the input variables $w = B'(x)p + q$ and $v = p + z$, we get $p = Q^{-1}(x)(diag(g(x))z - u - B(x)H(x)q)$. Consequently, the required control system is

$$\dot{x} = H(x)\left(I_n - B'(x)Q^{-1}(x)B(x)H(x)\right)q + H(x)B'(x)Q^{-1}(x)diag(g(x))z - H(x)B'(x)Q^{-1}(x)u$$

with inputs $q, z \in \Re^n, u \in \Re^n_+$. The computation of the feedback law for the above control system with Control Lyapunov Function $V(x) = \theta(x) - \theta(x^*)$, gives the dynamical system (2.6), where $F$ is defined by (2.4), (2.5). More specifically, we get:

$$\nabla V(x)\dot{x} = \nabla \theta(x)\left(H(x) - H(x)B'(x)Q^{-1}(x)B(x)H(x)\right)q$$
$$+ \nabla \theta(x)H(x)B'(x)Q^{-1}(x)diag(g(x))z - \nabla \theta(x)H(x)B'(x)Q^{-1}(x)u$$

The Control Lyapunov Function approach requires that each input must be selected so that each term appearing in the above equation takes negative values. The feedback laws

$$q = -R_1(x)\left(H(x) - H(x)B'(x)Q^{-1}(x)B(x)H(x)\right)(\nabla \theta(x))',$$
$$z = -(R_2(x)diag(g(x)) - diag(a(x)))Q^{-1}(x)B(x)H(x)(\nabla \theta(x))',$$
$$u = R_3(x)\left(Q^{-1}(x)B(x)H(x)(\nabla \theta(x))'\right)^+$$

where $R_1(x) \in \Re^{n \times n}$ is an arbitrary $C^1$, symmetric and positive definite matrix, $R_2(x) \in \Re^{k \times k}$ is an arbitrary $C^1$, symmetric and positive semidefinite matrix, $a_i(x)$, $b_i(x)$, $c_i(x)$ ($i = 1,...,k$) are arbitrary $C^1$ non-negative functions with $b_i(x) + c_i(x) > 0$ for all $i = 1,...,k$, $x \in S$ and at least one of the matrices $R_2(x) \in \Re^{k \times k}$, $diag(a(x)) \in \Re^{k \times k}$ being positive definite, where $a(x) := (a_1(x),...,a_k(x))' \in \Re^k$ and $R_3(x) \in \Re^{k \times k}$ is defined by (2.5), give us the vector field $F(x)$ defined by (2.4), (2.5).

**Remark 2.7:** If there are no equality constraints (i.e., $h(x) \equiv 0$) then the proof of Theorem 2.3 shows that exactly the same results with that of Theorem 2.3 hold with $H(x) \equiv I_n$.

**Proof of Theorem 2.3:** We first notice that statements (a) and (d) are direct consequences of definitions (2.3), (2.4), (2.5) and Fact 1. We next prove statements (b) and (c).

We first notice that definitions (2.5) and the fact that $g(x) \leq 0$ imply that the following equality holds for all $x \in S$:

$$\begin{aligned}
\nabla \theta(x) F(x) = &-\xi(x) R_1(x)(\xi(x))' - (diag(g(x))v(x))' R_2(x) diag(g(x))v(x) \\
&- \sum_{j=1}^{k} a_j(x) |g_j(x)| v_j^2(x) \\
&- \sum_{j=1}^{k} b_j(x) (\max(0, v_j(x)))^2 - \sum_{j=1}^{k} c_j(x) (\max(0, v_j(x)))^{2p_j+2} \\
\xi(x) = &\nabla \theta(x) \left[ H(x) - H(x) B'(x) Q^{-1}(x) B(x) H(x) \right]
\end{aligned} \quad (2.7)$$

where $\xi(x) = \nabla \theta(x) \left[ H(x) - H(x) B'(x) Q^{-1}(x) B(x) H(x) \right]$. It is clear that equation (2.7) shows that $\nabla \theta(x) F(x) \leq 0$ for all $x \in S$. We next investigate the nature of points $x \in S$ for which $\nabla \theta(x) F(x) = 0$. Equation (2.7) and the facts that $R_1(x) \in \Re^{n \times n}$ is a positive definite matrix, $R_2(x) \in \Re^{k \times k}$ is a positive semidefinite matrix, $a_i(x)$, $b_i(x)$, $c_i(x)$ ($i = 1,...,k$) are non-negative functions with $b_i(x) + c_i(x) > 0$ for all $i = 1,...,k$, $x \in S$ and at least one of the matrices $R_2(x) \in \Re^{k \times k}$, $diag(a(x)) \in \Re^{k \times k}$ is positive definite, where $a(x) := (a_1(x),...,a_k(x))' \in \Re^k$ and $R_3(x) \in \Re^{k \times k}$ is defined by (2.5), show that $\nabla \theta(x) F(x) = 0$ is equivalent to the following equations:

$$(v(x))^+ = 0 \quad (2.8)$$

$$g_j(x) v_j(x) = 0, \quad j = 1,...,k \quad (2.9)$$

$$\nabla \theta(x) H(x) \left[ I_n - B'(x) Q^{-1}(x) B(x) \right] H(x) = 0 \quad (2.10)$$

We define:

$$\mu = -v(x) \quad (2.11)$$

Definition (2.11) in conjunction with (2.8) and (2.9) implies that

$$\mu \geq 0 \text{ and } \mu' g(x) = 0 \quad (2.12)$$

Using (2.10) and the identity $H^2(x) = H(x)$, we obtain:

$$H(x) \left[ I_n - B'(x) Q^{-1}(x) B(x) H(x) \right] (\nabla \theta(x))' = 0 \quad (2.13)$$

Definitions (2.5), (2.11) in conjunction with (2.13) imply that:

$$H(x) (\nabla \theta(x) + \mu' B(x))' = 0 \quad (2.14)$$

Equation (2.14) in conjunction with Fact 3 and (2.12) implies that the conditions (1.3) hold. Therefore, $x \in \Phi$.

Thus, we have proved the implication: $\nabla \theta(x) F(x) = 0 \Rightarrow x \in \Phi$.

Consequently, we have proved statement (b) and one of the implications of statement (c) (namely, the implication $F(x) = 0 \Rightarrow x \in \Phi$).

We will prove now the implication $x \in \Phi \Rightarrow F(x) = 0$. Suppose that $x \in \Phi$. Then there exist $\lambda_i \in \Re$ ($i = 1,...,m$) and $\mu_j \geq 0$ ($j = 1,...,k$) such that conditions (1.3) hold, or in vector form:

$$\begin{aligned}(\nabla \theta(x))' + A'(x)\lambda + B'(x)\mu &= 0 \\ \mu' g(x) &= 0\end{aligned} \quad (2.15).$$

It follows from (2.15), Fact 1 and definitions (2.5) that

$$\begin{aligned}v(x) &= Q^{-1}(x) B(x) H(x) (\nabla \theta(x))' \\ &= -Q^{-1}(x) B(x) H(x) (A'(x)\lambda + B'(x)\mu) \\ &= -Q^{-1}(x) B(x) H(x) B'(x) \mu\end{aligned}$$

Using definition (2.3) and the above equality we obtain $v(x) = -\mu - Q^{-1}(x) diag(g(x))\mu$. However, the facts that $g(x) \leq 0$, $\mu \geq 0$ and $\mu' g(x) = 0$ imply that $diag(g(x))\mu = 0$. Consequently, it follows that $v(x) = -\mu$ and that (2.8) holds. Using (2.15), definition (2.4) and the facts that $v(x) = -\mu$, $diag(g(x))\mu = 0$, $(v(x))^+ = 0$, we obtain:

$$F(x) = -[H(x) - P'(x)Q(x)P(x)] R_1(x) [H(x) - P'(x)Q(x)P(x)](\nabla \theta(x))'$$

Using definitions (2.3), (2.5), Fact 1, (2.15), the above equality and the fact that $diag(g(x))\mu = 0$, we get:

$$\begin{aligned}F(x) &= H(x)[I_n - B'(x)Q^{-1}(x)B(x)]H(x)R_1(x)H(x)[I_n - B'(x)Q^{-1}(x)B(x)]H(x)(A'(x)\lambda + B'(x)\mu) \\ &= H(x)[I_n - B'(x)Q^{-1}(x)B(x)]H(x)R_1(x)H(x)[I_n - B'(x)Q^{-1}(x)B(x)]H(x)B'(x)\mu \\ &= H(x)[I_n - B'(x)Q^{-1}(x)B(x)]H(x)R_1(x)H(x)B'(x)[I_k - Q^{-1}(x)B(x)H(x)B'(x)]\mu \\ &= H(x)[I_n - B'(x)Q^{-1}(x)B(x)]H(x)R_1(x)H(x)B'(x)[I_k - Q^{-1}(x)(Q(x) + diag(g(x)))]\mu \\ &= -H(x)[I_n - B'(x)Q^{-1}(x)B(x)]H(x)R_1(x)H(x)B'(x)Q^{-1}(x)diag(g(x))\mu = 0\end{aligned}$$

We next turn to the proof of properties (1) and (2).

Local existence of the solution of the initial value problem (2.6) with $x(0) = x_0$ is a direct consequence of properties (a), (d) and the Nagumo theorem (page 27 in [3]). Global existence of the solution of the initial value problem (2.6) with $x(0) = x_0$ follows from Theorem 1.2.3 (page 27) in [3], assumption (H1) and the fact that $\theta(x(t))$ is non-increasing (a direct consequence of property (b)). In fact, assumption (H1) in conjunction with the fact that $\theta(x(t))$ is non-increasing shows that $\{x(t): t \geq 0\}$ is bounded.

As in the case of dynamical systems on $\Re^n$, it follows that $\omega(x_0)$ is a compact, positively invariant set for system (2.6) (see [24]). The fact that $\theta(x(t))$ is non-increasing implies that $\lim_{t \to +\infty} \theta(x(t)) = l = \inf\{\theta(x(t)): t \geq 0\}$, which shows that $\omega(x_0) \subseteq \{x \in S: \theta(x) = l\}$. We next show that $\omega(x_0) \subseteq \Phi$. The inequality

$$\int_0^t \gamma(x(s))ds \leq \theta(x_0) - l, \text{ for all } t \geq 0 \quad (2.16)$$

where $\gamma(x) := \nabla \theta(x) F(x)$ is a direct consequence of (2.7) and the definition $l := \inf\{\theta(x(t)): t \geq 0\}$. Notice that the mapping $\Re_+ \ni s \to \gamma(x(s))$ is uniformly continuous since $\{\dot{x}(t) = F(x(t)): t \geq 0\}$ is bounded (a consequence of the fact that $\{x(t): t \geq 0\}$ is bounded) and since the mapping $\Re^n \ni x \to \gamma(x)$ is locally Lipschitz. Using (2.16) and applying Barbalat's lemma (see [18]), we conclude that $\lim_{s \to +\infty} \gamma(x(s)) = 0$. The validity of the implication $\nabla \theta(x) F(x) = \gamma(x) = 0 \Rightarrow x \in \Phi$ implies that $\omega(x_0) \subseteq \Phi$.

Finally, the fact that every strict local solution $x^* \in S$ of the NLP problem described by (1.1) and (1.2) follows from property (b) and the consideration of the Lyapunov function $V(x) = \theta(x) - \theta(x^*)$. The proof is complete. ◁

## 3. Numerical solutions of NLP problems

As remarked in the Introduction and in [5], each system of differential equations that solves a NLP problem when combined with a numerical scheme for solving Ordinary Differential Equations (ODEs) provides a numerical scheme for solving the NLP problem. However, when we try to apply a numerical scheme for the solution of (2.6) then we face the problem that the dynamical system (2.6) is not defined on $\Re^n$ but on the closed set $S \subseteq \Re^n$.

In the literature, projection schemes have been proposed (see [13,14]). The projection schemes preserve the order of the applied numerical scheme (see [13,14]) even if the projection on the closed set $S \subseteq \Re^n$ is not exact. However, the application of a Runge-Kutta numerical scheme for (2.6) and its (approximate) projection on the closed set $S \subseteq \Re^n$ means that the solution of a NLP problem is required at each time step. The corresponding NLP problem may be as difficult as the initial one, which means that this approach is not easily applicable (with the exception of cases where the projection is easy, see [12]).

The key idea presented in this work is that the selection of the applied time step can be used for solving the above problems. First we focus on the case without equality constraints.

The following theorem is the main result of this section, which guarantees global convergence of the above algorithm.

**Theorem 3.1:** *Suppose that assumptions (H1), (H2) hold for the NLP problem described by (1.1) and (1.2) with $h(x) \equiv 0$. Let $R_1(x) \in \Re^{n \times n}$ be an arbitrary $C^1$, symmetric and positive definite matrix, $R_2(x) \in \Re^{k \times k}$ be an arbitrary $C^1$, symmetric and positive semidefinite matrix, $a_i(x)$, $b_i(x)$, $c_i(x)$ ($i=1,...,k$) be arbitrary $C^1$ non-negative functions with $b_i(x) + c_i(x) > 0$ for all $i=1,...,k$, $x \in S$ and at least one of the matrices $R_2(x) \in \Re^{k \times k}$, $diag(a(x)) \in \Re^{k \times k}$ being positive definite, where $a(x) := (a_1(x),...,a_k(x))' \in \Re^k$. Let $F(x) \in \Re^n$ be the vector field defined by (2.4), (2.5) with $H(x) \equiv I_n$. Furthermore, we assume that:*

*Consider the following algorithm:*

**Algorithm:** Given constants $r > 0$, $\varepsilon \in (0,r)$, $\lambda \in (0,1)$ and an initial point $x_0 \in S$, we follow the steps for $i=0,1,...$

→ Step i: Calculate $F(x_i)$ using (2.4), (2.5). If $|F(x_i)| = 0$ then $x_{i+1} = x_i$. If $|F(x_i)| > 0$ then set $s^{(0)} = r$ and $p = 0$. Moreover, let $I(x_i) \subseteq \{1,...,k\}$ be the set of all indices $j \in \{1,...,k\}$ with $\max_{0 \leq s \leq \varepsilon}(g_j(x_i + s F(x_i))) > -\varepsilon$.

  → Step p: Calculate $x_i^{(p)} = x_i + s^{(p)} F(x_i)$.

  → Solve $\min\left\{\left|y - x_i^{(p)}\right|^2 : \max_{j \in I(x_i)}(g_j(y)) \leq 0\right\}$, for the case $I(x_i) \neq \varnothing$ or set $y = x_i^{(p)}$, for the case $I(x_i) = \varnothing$.

    → If $y \in S$ and $\theta(y) \leq \theta(x_i) + \lambda s^{(p)} \nabla \theta(x_i) F(x_i)$ then set $x_{i+1} = y$, $i = i+1$ and go to Step i.

    → If $y \notin S$ or $\theta(y) > \theta(x_i) + \lambda s^{(p)} \nabla \theta(x_i) F(x_i)$ then set $s^{(p+1)} = \frac{1}{2} s^{(p)}$, $p = p+1$ and go to Step p.

*Then every accumulation point $x^*$ of the sequence $x_i$ produced by the above algorithm satisfies $x^* \in \Phi$.*

**Remark 3.2:** It is clear that the algorithm presented in Theorem 3.1 exploits the time step used for the estimate provided by the explicit Euler scheme $x_i^{(p)} = x_i + s^{(p)} F(x_i)$. The constant $r > 0$ is the maximum allowable time step. In most iterations, the algorithm will not require the solution of an NLP problem, provided that $\varepsilon > 0$ is sufficiently small. The fact that in most cases the Euler scheme is sufficient is explained by statement (d) of Theorem 2.3: for all $j = 1,...,k$ it holds that

$$\nabla g_j(x) F(x) = g_j(x) \omega_j(x) - \left(b_j(x) + c_j(x)(\max(0, v_j(x)))^{2 p_j}\right) \max(0, v_j(x)) \tag{3.1}$$

where $\omega(x) := (\omega_1(x),...,\omega_k(x))'$ is given by

$$\begin{aligned}\omega(x) := &\, Q^{-1}(x) B(x) H(x) R_1(x) [H(x) - P'(x) Q(x) P(x)] (\nabla \theta(x))' \\ &- [I_k + Q^{-1}(x) \operatorname{diag}(g(x))] (R_2(x) \operatorname{diag}(g(x)) - \operatorname{diag}(a(x))) v(x) - Q^{-1}(x) R_3(x) (v(x))^+ \end{aligned} \tag{3.2}$$

Using the estimate

$$g_j(x + sF(x)) \leq g_j(x) + s \nabla g_j(x) F(x) + \frac{1}{2} s^2 K_j(x) \tag{3.3}$$

where $s \in [0, r]$, $K_j(x) := \max\{F'(x) \nabla^2 g_j(x + \mu F(x)) F(x) : \mu \in [0, r]\}$, it follows that $g_j(x + sF(x)) \leq -\varepsilon$ provided that

$$\varepsilon + g_j(x) + s g_j(x) \omega_j(x) - s \left(b_j(x) + c_j(x)(\max(0, v_j(x)))^{2 p_j}\right) \max(0, v_j(x)) + \frac{1}{2} s^2 K_j(x) \leq 0$$

The above inequality is satisfied for all $s \in [0, \varepsilon]$ in many cases provided that $\varepsilon > 0$ is sufficiently small. This explains the additional fact that that the NLP problem $\min\left\{ \left|y - x_i^{(p)}\right|^2 : \max_{j \in I(x_i)} (g_j(y)) \leq 0 \right\}$ is much simpler than the initial NLP problem: the index set $I(x_i) \subseteq \{1,...,k\}$ is expected to be a set with small cardinal number.

Finally, as remarked in [12], the solution of the NLP problem $\min\left\{ \left|y - x_i^{(p)}\right|^2 : \max_{j \in I(x_i)} (g_j(y)) \leq 0 \right\}$ need not be exact: it suffices to assume that to find any point $y \in \Re^n$ with $\max_{j \in I(x_i)} (g_j(y)) \leq 0$ and $\left|y - x_i^{(p)}\right| \leq C \left|y^* - x_i^{(p)}\right|$, where $C \geq 1$ is a constant and $y^* \in \Re^n$ is any of the global solutions of the NLP problem $\min\left\{ \left|y - x_i^{(p)}\right|^2 : \max_{j \in I(x_i)} (g_j(y)) \leq 0 \right\}$.

**Proof:** Define the sets:

$$\tilde{\Phi} := \{ z \in \Phi, \theta(z) \leq \theta(x_0) \}, \quad \tilde{S} := \{ z \in S, \theta(z) \leq \theta(x_0) \} \tag{3.4}$$

Notice that the set $\tilde{\Phi}$ is non-empty and compact (by virtue of property (b) of Theorem 2.3 it follows that the set $\tilde{\Phi}$ coincides with the set $\{x \in S : |\nabla \theta(x) F(x)| = 0, \theta(x) \leq \theta(x_0)\}$ for which assumption (H1) implies that it is bounded; notice that the set $\tilde{\Phi}$ contains the global solution of the NLP problem described by (1.1) and (1.2)). Moreover, assumption (H1) guarantees that $\tilde{S}$ is non-empty and compact.

Let $d(x)$ be the distance of any point $x \in \Re^n$ from the set $\tilde{\Phi} := \{ y \in \Phi : \theta(y) \leq \theta(x_0) \}$, i.e.,

$$d(x) := \inf \{ |x - y| : y \in \tilde{\Phi} \} \tag{3.5}$$

Since the set $\tilde{\Phi}$ is non-empty and compact, it follows that the function $d(x)$ is well-defined, is globally Lipschitz (with unit Lipschitz constant) and satisfies $d(x) > 0$ for all $x \notin \tilde{\Phi}$.

**Claim:** *For every $\eta > 0$ there exists a constant $\delta_\eta > 0$ such that the following implication holds:*

"If $x \in \widetilde{S}$, $d(x) \geq \eta$ and $s \leq \delta_\eta$ then $y \in S$ and $\theta(y) \leq \theta(x) + \lambda s \nabla \theta(x) F(x)$,

where $y = x + sF(x)$ for the case $I(x) = \varnothing$ and

$y \in \Re^n$ is any global solution of $\min\left\{ |y - x - sF(x)|^2 : \max_{j \in I(x)}(g_j(y)) \leq 0 \right\}$ for the case $I(x) \neq \varnothing$" (3.6)

The proof of the claim can be found at the Appendix.

We notice that, by virtue of implication (3.6) and the fact that $\theta(x_i)$ is non-increasing, the algorithm is well-defined (i.e., for each iteration $i$ the variable $p$ assumes finite values). Let $s_i \in (0, r]$ ($i = 0, 1, ...$) be the applied $s^{(p)}$ for which $y \in S$ and $\theta(y) \leq \theta(x_i) + \lambda s^{(p)} \nabla \theta(x_i) F(x_i)$ for the case $x_i \notin \Phi$ and $s_i = 0$ for the case $x_i \in \Phi$. For every $i = 0, 1, ...$ it holds that:

$$\theta(x_{i+1}) \leq \theta(x_i) + \lambda s_i \nabla \theta(x_i) F(x_i) \quad (3.7)$$

Notice that assumption (H1) in conjunction with (3.5) guarantees that the sequence $x_i$ is bounded with $\theta(x_i) \leq \theta(x_0)$ for all $i = 0, 1, ...$. Moreover, implication (3.6) implies the following inequality for every $i = 0, 1, ...$ with $x_i \notin \Phi$:

$$s_i \geq \delta_\eta / 2, \text{ for all } \eta \in (0, d(x_i)] \quad (3.8)$$

where $\delta_\eta > 0$ is the constant involved in implication (3.6). If $x^*$ is one global solution of the NLP problem described by (1.1) and (1.2) then by applying (3.7) inductively, we conclude that the following inequality holds for $i = 1, 2, ...$

$$\lambda \sum_{l=0}^{i-1} s_l |\nabla \theta(x_l) F(x_l)| \leq \theta(x_0) - \theta(x^*) \quad (3.9)$$

The above inequality implies that $\lim_{i \to +\infty} (s_i |\nabla \theta(x_i) F(x_i)|) = 0$.

In order to prove that every accumulation point $x^*$ of the sequence $x_i$ produced by the above algorithm satisfies $x^* \in \Phi$, we will use a contradiction argument. Let a subsequence of the sequence $x_i$ which converges. We will use the same notation $x_i$ for the subsequence and let $x^*$ be the unique accumulation point of the subsequence $x_i$. We assume that $x^* \notin \Phi$. By continuity and using property (b) of Theorem 2.3, we have $\lim_{i \to +\infty} (|\nabla \theta(x_i) F(x_i)|) = |\nabla \theta(x^*) F(x^*)| > 0$. Since $\lim_{i \to +\infty} (s_i |\nabla \theta(x_i) F(x_i)|) = 0$, we are in a position to conclude that

$$\lim_{i \to +\infty} (s_i) = 0 \quad (3.10)$$

Since $x^* \notin \widetilde{\Phi}$, it follows that $\lim_{i \to +\infty}(d(x_i)) = d(x^*) > 0$. Therefore, there exists $\eta > 0$ such that $d(x_i) \geq \eta$ for sufficiently large $i$. Thus, inequality (3.8) gives $s_i \geq \delta_\eta / 2$, where $\delta_\eta > 0$ is the constant involved in implication (3.6). This contradicts (3.10).

The proof is complete. ◁

**Remark 3.3:** When equality constraints are present, then it should be noticed that Theorem 3.1 is still useful under the following assumption:

**(H3)** *There exist positive integers $n_1, n_2$ with $n_1 + n_2 = n$ and a smooth function $\phi: \Re^{n_1} \to \Re^{n_2}$ such that for every $\xi \in \Re^{n_1}$ it holds that $h(x) = 0$, where $x = (\xi, \phi(\xi))$.*

Indeed, under assumption (H3), we may define for all $\xi \in \Re^{n_1}$:

$$\tilde{\theta}(\xi) = \theta(x) \text{ with } x = (\xi, \phi(\xi)) \tag{3.11}$$

$$\tilde{g}_j(\xi) := g_j(x) \text{ with } x = (\xi, \phi(\xi)) \text{ for } j = 1,...,k \tag{3.12}$$

We can also define $\tilde{F}(\xi)$ to be the vector field that is made up from the first $n_1$ components of the vector field $F(x)$ defined by (2.4), (2.5) evaluated at $x = (\xi, \phi(\xi))$. Then we can apply Theorem 3.1 with $\tilde{\theta}, \tilde{g}_j (j=1,...,k), \tilde{F}$ in place of $\theta, g_j (j=1,...,k), F$.

**Remark 3.4:** The algorithm may be modified in a straightforward way for other higher order explicit Runge-Kutta numerical schemes. This is meaningful only when the vector field $F(x)$ has sufficient regularity. More specifically, the term $x_i^{(p)} = x_i + s^{(p)} F(x_i)$ may be replaced $x_i^{(p)} = N(s^{(p)}, x_i)$, for an appropriate mapping $N(s^{(p)}, x_i)$ which is characteristic of the Runge-Kutta scheme and the definition of the set $I(x_i) \subseteq \{1,...,k\}$ is modified to be the set of all indices $j \in \{1,...,k\}$ with $\max_{0 \le s \le \varepsilon} (g_j(N(s, x_i))) > -\varepsilon$. However, it should be noticed that for higher order explicit Runge-Kutta schemes, the vector field $F(x)$ must be computed for various points. Since $F(x)$ is defined only for a neighborhood of the set $S$, it may be needed to restrict the time step $s^{(p)}$ so that all points which are needed for the evaluation of $N(s^{(p)}, x_i)$ are in a neighborhood of the set $S$.

**Remark 3.5:** Using (3.1), (3.2) and (3.3) and assuming that for all $j = 1,...,k$, there exist positive, continuous functions $q_j : S \to (0, +\infty)$, $Q_j : S \to (0, +\infty)$ ($j = 1,...,k$), such that the following inequalities hold for all $j = 1,...,k$ and $x \in S$:

$$K_j(x) \le -g_j(x) q_j(x) + Q_j(x) \left( b_j(x) + c_j(x) (\max(0, v_j(x)))^{2p_j} \right) \max(0, v_j(x)) \tag{3.13}$$

where $K_j(x) := \max\{F'(x) \nabla^2 g_j(x + \mu F(x)) F(x) : \mu \in [0, r]\}$ ($j = 1,...,k$), we can conclude that

$$g_j(x + sF(x)) \le 0, \text{ for all } j = 1,...,k \text{ and } x \in S \tag{3.14}$$

provided that

$$s \le \min_{j=1,...,k} \left( r, \frac{\omega_j(x) + \sqrt{\omega_j^2(x) + 2q_j(x)}}{q_j(x)}, \frac{2}{Q_j(x)} \right) \tag{3.15}$$

Define:

$$s_g(x) := \sup \left\{ s \in [0, r] : \max_{\substack{0 \le l \le s \\ j=1,...,k}} (g_j(x + lF(x))) \le 0 \right\}, \text{ for all } x \in S \tag{3.16}$$

and notice that (3.15) implies $s_g(x) \ge \min_{j=1,...,k} \left( r, \frac{\omega_j(x) + \sqrt{\omega_j^2(x) + 2q_j(x)}}{q_j(x)}, \frac{2}{Q_j(x)} \right)$ for all $x \in S$. Using the analogue of inequality (3.3) for $\theta(x)$, i.e., the inequality

$$\theta(x + sF(x)) \le \theta(x) + s \nabla \theta(x) F(x) + s^2 K_\theta(x)/2 \tag{3.17}$$

where $s \in [0, r]$, $K_\theta(x) := \max\{F'(x) \nabla^2 \theta(x + \mu F(x)) F(x) : \mu \in [0, r]\}$, we can conclude that the best possible choice for the time step $s \in [0, r]$ is given by:

$$s = \min \left( s_g(x), \frac{|\nabla \theta(x) F(x)|}{K_\theta(x)} \right), \text{ for the case } K_\theta(x) > 0 \text{ and } s = s_g(x), \text{ for the case } K_\theta(x) \le 0 \tag{3.18}$$

We notice that inequalities (3.13) hold automatically for arbitrary positive, continuous functions $q_j : S \to (0,+\infty)$, $Q_j : S \to (0,+\infty)$ ($j=1,...,k$), for the case where all functions $g_j(x)$ ($j=1,...,k$) are linear.

Therefore, if inequality (3.13) holds then we can simply compute the sequence $x_{i+1} = M(x_i)$, where $M(x) = x + sF(x)$ and $s \in [0,r]$ is given by (3.18). We notice that the implementation of an approximation of the numerical scheme $x_{i+1} = M(x_i)$ does not necessarily requires knowledge of the functions $K_j(x)$, $q_j : S \to (0,+\infty)$, $Q_j : S \to (0,+\infty)$ ($j=1,...,k$) and $K_\theta(x)$: evaluating $g(x+sF(x))$ and $\theta(x+sF(x))$ for certain $s \in [0,r]$ can give us estimates of $K_j(x)$ and $K_\theta(x)$ satisfying (3.3) and (3.17). Using these estimates we can estimate $s_g(x)$. Consequently, the algorithm is implemented as follows:

**Algorithm:** Given constants $r > 0$, $\varepsilon > 0$, $\lambda \in (0,1/2]$ and an initial point $x_0 \in S$, we follow the steps for $i = 0,1,...$

→ Step i: Calculate $F(x_i)$ using (2.4), (2.5). If $|F(x_i)| = 0$ then $x_{i+1} = x_i$. If $|F(x_i)| > 0$ then $K_j^{(0)} = \max\left(\varepsilon, 2r^{-2}\left(g_j(x_i + rF(x_i)) - g_j(x_i) - r\nabla g_j(x_i)F(x_i)\right)\right)$ for $j=1,...,k$, $K_\theta^{(0)} = \max\left(\varepsilon, 2r^{-2}\left(\theta(x+rF(x)) - \theta(x) - r\nabla\theta(x)F(x)\right)\right)$ and set $p = 0$.

→ Step p: Compute $s_j^{(p)} = \min\left(r, \dfrac{-\nabla g_j(x_i)F(x_i) + \sqrt{|\nabla g_j(x_i)F(x_i)|^2 - 2K_j^{(p)}g_j(x_i)}}{K_j^{(p)}}\right)$ for $j=1,...,k$ and

$s^{(p)} = \min\limits_{j=1,...,k}\left(s_j^{(p)}, \dfrac{|\nabla\theta(x_i)F(x_i)|}{K_\theta^{(p)}}\right)$. Calculate $x_i^{(p)} = x_i + s^{(p)}F(x_i)$.

→ If $x_i^{(p)} \in S$ and $\theta(x_i^{(p)}) \leq \theta(x_i) + \lambda s^{(p)}\nabla\theta(x_i)F(x_i)$ then set $x_{i+1} = x_i^{(p)}$, $i = i+1$ and go to Step i.

→ If $x_i^{(p)} \notin S$ or $\theta(x_i^{(p)}) > \theta(x_i) + \lambda s^{(p)}\nabla\theta(x_i)F(x_i)$ then set $K_j^{(p+1)} = K_j^{(p)} + \varepsilon$ ($j = 1,...,k$), $K_\theta^{(p+1)} = K_\theta^{(p)} + \varepsilon$, $p = p+1$ and go to Step p.

Using exactly the same procedure with that of the proof of Theorem 3.1, we can conclude that every accumulation point $x^*$ of the sequence $x_i$ produced by the above algorithm satisfies $x^* \in \Phi$, provided that assumptions (H1), (H2) and (3.13) hold. However, numerical experiments show that the algorithm can converge even when assumption (3.13) does not hold.

## 4. Examples

In order to demonstrate the performance of the proposed algorithms we have used two examples from [28]. The first example is dealing with the solution of the problem:

$$\min \theta(x) = x_1^2 + 2x_2^2 + x_1 x_2 - 6x_1 - 2x_2 - 12x_3$$
$$s.t.$$
$$h(x) = x_1 + x_2 + x_3 - 2 = 0$$
$$g(x) = \begin{bmatrix} -x_1 + 2x_2 - 3 \\ -x_1 \\ -x_2 \\ -x_3 \end{bmatrix} \leq 0$$
(4.1)

It can be shown that assumptions (H1), (H2) hold for this problem. Moreover, assumption (H3) holds with $\xi = (x_1, x_2)' \in \Re^2$ and $\phi(\xi) = 2 - x_1 - x_2$. Since, the inequality constraints are linear, we are in a position to use the algorithm of Remark 3.5. We have used the algorithm of Remark 3.5 with $R_1(x) = \sigma I_3$, where $\sigma > 0$,

$R_2(x) = 0 \in \Re^{4\times 4}$, $a_i(x) \equiv 1$, $b_i(x) \equiv 1$, $c_i(x) \equiv 0$ ($i=1,...,4$), $r=1$, $\lambda = 0.1$ and $\varepsilon = 10^{-6}$.

It was found that for all initial points in the feasible set and for every $\sigma \in [0.01, 200]$ the algorithm converges at the point $(x_1, x_2, x_3) = (0,0,2)$ in no more than 3 iterations. In this case, the convergence of the algorithm of Remark 3.5 is very fast.

Figure 1 shows the projection of the phase diagram on the $x_1 - x_2$ plane for the dynamical system (2.6), where $F$ is defined by (2.4), (2.5), for problem (4.1) with $R_1(x) = \sigma I_3$, $\sigma = 2$, $R_2(x) = 0 \in \Re^{4\times 4}$, $a_i(x) \equiv 1$, $b_i(x) \equiv 1$, $c_i(x) \equiv 0$ ($i=1,...,4$). Figure 1 was created by solving numerically system (2.6) with the explicit Euler method and time step 0.01.

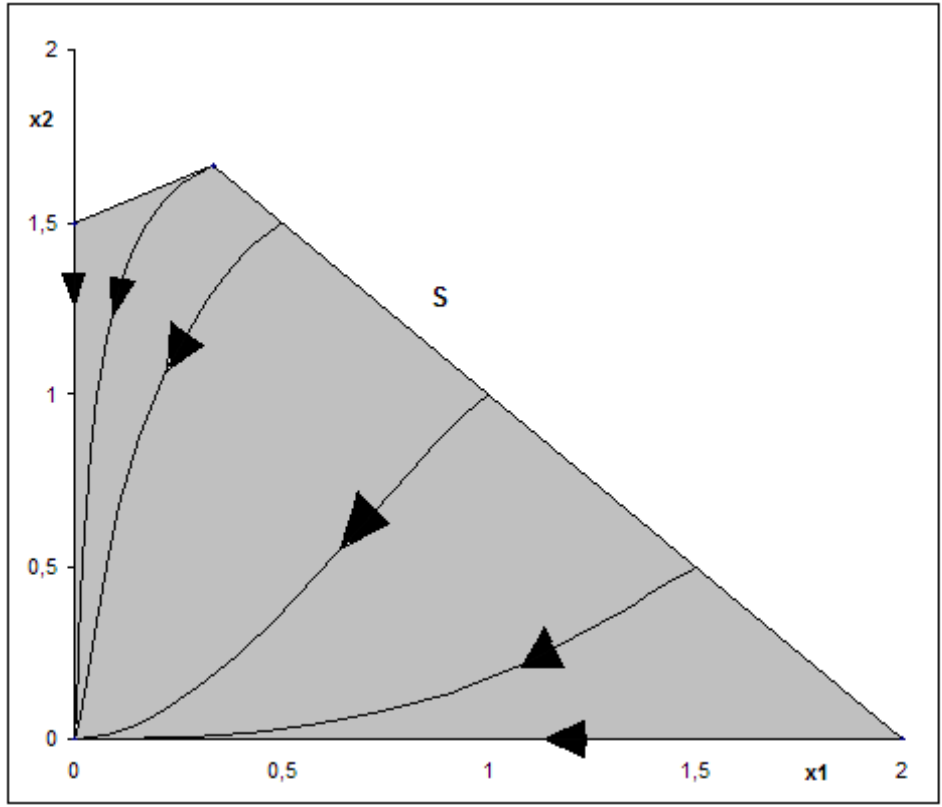

Figure 1: The projection of the phase diagram on the $x_1 - x_2$ plane for the dynamical system (2.6), where $F$ is defined by (2.4), (2.5), for problem (4.1) with $R_1(x) = \sigma I_3$, $\sigma = 2$, $R_2(x) = 0 \in \Re^{4\times 4}$, $a_i(x) \equiv 1$, $b_i(x) \equiv 1$, $c_i(x) \equiv 0$ ($i=1,...,4$).

The second example is dealing with the Rosen–Suzuki problem:

$$\min \theta(x) = x_1^2 + x_2^2 + 2x_3^2 + x_4^2 - 5x_1 - 5x_2 - 21x_3 + 7x_4$$
s.t.
$$h(x) = 2x_1^2 + x_2^2 + x_3^2 + 2x_1 - x_2 - x_4 - 5 = 0 \qquad (4.2)$$
$$g(x) = \begin{bmatrix} x_1^2 + x_2^2 + x_3^2 + x_4^2 + x_1 - x_2 + x_3 - x_4 - 8 \\ x_1^2 + 2x_2^2 + x_3^2 + 2x_4^2 - x_1 - x_4 - 10 \end{bmatrix} \leq 0$$

It can be shown that assumptions (H1), (H2) hold for this problem. Moreover, assumption (H3) holds with $\xi = (x_1, x_2, x_3)' \in \Re^3$ and $\phi(\xi) = 2x_1^2 + x_2^2 + x_3^2 + 2x_1 - x_2 - 5$. The vector field $F(x)$ defined by (2.4), (2.5) is constructed with $R_1(x) = \sigma I_4$, where $\sigma > 0$, $R_2(x) = 0 \in \Re^{2\times 2}$, $a_i(x) \equiv 1$, $b_i(x) \equiv 1$, $c_i(x) \equiv 0$ ($i=1,2$).

This is a problem with nonlinear inequality constraints. Therefore, we cannot assume the validity of (3.13). Indeed, there are points in the feasible set with $g_1(x) = 0$, $\max(0, v_1(x)) = 0$ and for which $\tilde{g}_1(\xi + s\tilde{F}(\xi)) > 0$ for $s > 0$, where $\tilde{g}_1$ is defined by (3.12) and $\tilde{F}(\xi)$ is the vector field that is made up from the first 3 components of the vector field $F(x)$ evaluated at $x = (\xi, \phi(\xi))$. Such a point is $x = (-1, -1, 2, 1)' \in \Re^4$ and it is clear that we cannot apply the algorithm of Remark 3.5 at any one of these points. However, the algorithm of Remark 3.5 may be applied with other initial points: for example, if the algorithm of Remark 3.5 is applied with $\sigma = 0.2$, $r = 1$, $\lambda = 0.1$ and $\varepsilon = 10^{-6}$ to the initial point $x_0 = (-0.9, -1, 2, 0.82)' \in \Re^4$ (which is close to the "problematic point" $x = (-1, -1, 2, 1)' \in \Re^4$) then the produced sequence reaches the neighborhood $N = \{|x - x^*| \leq 10^{-5}\}$, with $x^* = (0, 1, 2, -1)' \in \Re^4$ in 33 iterations. It was also found that different values of $\sigma > 0$ and $r > 0$ affect the convergence properties of the algorithm. For example, lower values than 1 for $r > 0$ and higher values than 0.5 for $\sigma > 0$ require more iterations for convergence. The algorithm of Remark 3.5 performs well from almost all points of the feasible set: for example, if the algorithm of Remark 3.5 is applied with $\sigma = 0.2$, $r = 1$, $\lambda = 0.1$ and $\varepsilon = 10^{-6}$ to the initial point $x_0 = (-1, -1, -2, 1)' \in \Re^4$ then the produced sequence reaches the neighborhood $N = \{|x - x^*| \leq 10^{-5}\}$, with $x^* = (0, 1, 2, -1)' \in \Re^4$ in 47 iterations.

For the initial point $x_0 = (-1, -1, 2, 1)' \in \Re^4$ we can apply the algorithm of Theorem 3.1. If the algorithm of Theorem 3.1 is applied with $\sigma = 0.2$, $r = 0.5$, $\lambda = 0.1$ and $\varepsilon = 10^{-6}$ to the initial point $x_0 = (-1, -1, 2, 1)' \in \Re^4$ then the produced sequence reaches the neighborhood $N = \{|x - x^*| \leq 10^{-5}\}$, with $x^* = (0, 1, 2, -1)' \in \Re^4$ in 39 iterations.

In general, the convergence of the proposed algorithms is linear. For superlinear convergence, we can either use different selections for $R_1(x) \in \Re^{4 \times 4}$, $R_2(x) \in \Re^{2 \times 2}$, $a_i(x)$, $b_i(x)$, $c_i(x)$ ($i = 1, 2$) or use a different algorithm once we are close to the set $\Phi$. The quantity $|F(x)|$ can be used in order to signal the approach of a neighborhood of $\Phi$.

## 5. Conclusions

In this work we have showed that given a nonlinear programming problem, it is possible, under mild assumptions, to construct a family of dynamical systems defined on the feasible set of the given problem, so that: (a) the equilibrium points are the unknown critical points of the problem, (b) each dynamical system admits the objective function of the problem as a Lyapunov function, and (c) explicit formulae are available without involving the unknown critical points of the problem. The construction of the family of dynamical systems is based on the Control Lyapunov Function methodology, which is used in mathematical control theory for the construction of stabilizing feedback.

The knowledge of a dynamical system with the previously mentioned properties allows the construction of algorithms which guarantee global convergence to the set of the critical points. However, we make no claim about the effectiveness of the proposed algorithms. The topic of the numerical solution of NLPs is a mature topic and it is clear that other algorithms have much better characteristics than the algorithms proposed in this paper. However, the theory used for the construction of the algorithm is different from other existing algorithms. The algorithms contained in this work are derived by using concepts of dynamical systems theory and mathematical control theory.

The obtained results have nothing to do with extremum seeking (see [10,19]), but may open the way of using different extremum seeking control schemes in the future for constrained problems. Finally, it may be beneficial to compare the algorithm with other global algorithms (see [22] and references therein): this is a future research topic.

**Acknowledgements:** The author would like to thank Prof. Lars Grüne for his valuable comments and suggestions. The contribution of Prof. Lars Grüne to this work is major.

# Appendix

**Proof of Lemma 2.2:** First notice that by virtue of Fact 2, the following equality holds for all $\xi = (\xi_1, ..., \xi_k)' \in \Re^k$:

$$\xi' Q(x) \xi = |H(x) B'(x) \xi|^2 + \sum_{j=1}^{k} |g_j(x)| \xi_j^2 \qquad (A.1)$$

Equality (A.1) implies that $Q(x) \in \Re^{k \times k}$ is positive semidefinite. Suppose that $Q(x) \in \Re^{k \times k}$ is not positive definite.

Then there exists a non-zero $\xi = (\xi_1,...,\xi_k)' \in \Re^k$ with $\xi'Q(x)\xi = 0$. Consequently, equality (A.1) shows that we must have $H(x)B'(x)\xi = 0$ and $\xi_j = 0$ for all $j = 1,...,k$ with $g_j(x) < 0$. Fact 3 implies that there exists $\lambda \in \Re^m$ such that $B'(x)\xi = A'(x)\lambda$. The previous equality implies that

$$\sum_{j=1}^{k} \xi_j \nabla g_j(x) - \sum_{i=1}^{m} \lambda_i \nabla h_i(x) = 0 \tag{A.2}$$

Since $\xi_j = 0$ for all $j = 1,...,k$ with $g_j(x) < 0$ and since $\xi = (\xi_1,...,\xi_k)' \in \Re^k$ is non-zero, we conclude from (A.2) that assumption (H2) is violated.

The proof is complete. ◁

**Proof of the Claim:** Let $\eta > 0$ be arbitrary. We distinguish two cases.

Case 1: The set $\{x \in \widetilde{S} : d(x) \geq \eta\}$ is empty, where $\widetilde{S} \subseteq S$ is defined by (3.4) and $d(x)$ is defined by (3.5). In this case, implication (3.6) holds trivially with arbitrary $\delta_\eta > 0$.

Case 2: The set $\{x \in \widetilde{S} : d(x) \geq \eta\}$ is non-empty.

Continuity of the distance function $d(x)$ and compactness of $\widetilde{S} \subseteq S$ implies that the set $\{x \in \widetilde{S} : d(x) \geq \eta\}$ is compact. Statements (b) and (c) of Theorem 2.3 guarantee that the quantity

$$\rho := \min\left\{\frac{|\nabla \theta(x) F(x)|}{|F(x)|(|\nabla \theta(x)| + |F(x)|)} : x \in \widetilde{S}, d(x) \geq \eta\right\} \tag{A.3}$$

is well-defined and is positive.

Let $x \in \widetilde{S}$ with $d(x) \geq \eta$ be an arbitrary point. We denote by $z(t)$ the unique solution of the initial value problem $\dot{z} = F(z)$ with $z(0) = x$. We also notice that the vector field $F$ as defined by (2.4), (2.5) is locally Lipschitz on a neighborhood of $S$. By virtue of compactness of $\widetilde{S} \subseteq S$, we are in a position to assume the existence of a constant $L \geq 0$ that satisfies:

$$|F(y) - F(x)| \leq L|y - x|, \text{ for all } x, y \in \widetilde{S} \tag{A.4}$$

Inequality (A.4), the fact that $z(s)$ belongs to the compact set of all $z \in S$ with $\theta(z) \leq \theta(x)$ and standard arguments show that the following inequality holds for all $s \geq 0$:

$$|z(s) - x - sF(x)| \leq Le^{Ls}\frac{s^2}{2}|F(x)| \tag{A.5}$$

Next we notice that the problem $\min\left\{|y - x - sF(x)|^2 : \max_{j \in I(x)}(g_j(y)) \leq 0\right\}$ with $I(x) \neq \emptyset$ admits at least one solution (since the mapping $y \to |y - x - sF(x)|^2$ is radially unbounded). Any solution $y \in \Re^n$ of the problem $\min\left\{|y - x - sF(x)|^2 : \max_{j \in I(x)}(g_j(y)) \leq 0\right\}$ with $I(x) \neq \emptyset$ satisfies for all $s \geq 0$:

$$\begin{aligned} |y - x - sF(x)| &\leq |z(s) - x - sF(x)| \\ |y - x - sF(x)| &\leq s|F(x)| \end{aligned} \tag{A.6}$$

Notice that the above inequalities hold trivially for the case $I(x) = \emptyset$ and $y = x + sF(x)$. Define:

$$q := \max\left\{\sum_{j=1}^{k}|\nabla g_j(z)| : z \in \Re^n, |z| \leq \beta\right\}, \text{ where } \beta := \max\left\{|x| + 2r|F(x)| : x \in \widetilde{S}, d(x) \geq \eta\right\} \quad (A.7)$$

$$Q := \frac{Le^{Lr}}{2} + 2\max\left\{|\nabla^2 \theta(z)| : z \in \Re^n, |z| \leq \beta\right\} \quad (A.8)$$

We will show next that implication (3.6) holds with $\delta_\eta > 0$ defined by:

$$\delta_\eta := \min\left(\varepsilon, \left(\frac{2\varepsilon}{1 + qL\gamma e^{Lr}}\right)^{1/2}, \frac{(1-\lambda)\rho}{1+Q}\right), \text{ where } \gamma := \max\left\{|F(x)| : x \in \widetilde{S}, d(x) \geq \eta\right\} \quad (A.9)$$

First, we show the implication:

$$\text{"If } x \in \widetilde{S}, \ d(x) \geq \eta \text{ and } s \leq \delta_\eta \text{ then } y \in S\text{"} \quad (A.10)$$

It suffices to show that $g_j(y) \leq 0$ for all $j \notin I(x)$. Notice that by virtue of the definition of the set $I(x) \subseteq \{1,...,k\}$ it follows that $g_j(x + sF(x)) \leq -\varepsilon$, for all $s \in [0, \varepsilon]$ and $j \notin I(x)$. Using (A.6) we obtain for all $s \in [0, r]$:

$$g_j(y) \leq g_j(x + sF(x)) + |y - x - sF(x)|\max\left\{|\nabla g_j(z)| : |z - x| \leq 2r|F(x)|\right\} \quad (3.11)$$

Since $g_j(x + sF(x)) \leq -\varepsilon$, for all $s \in [0, \varepsilon]$ and $j \notin I(x)$, we obtain from (A.5), (A.6), (A.11) and (A.7) for all $s \in [0, \varepsilon]$ and $j \notin I(x)$:

$$g_j(y) \leq -\varepsilon + qL|F(x)|e^{Lr}\frac{s^2}{2} \quad (A.12)$$

Inequality (A.12) in conjunction with definition (A.9) shows that $g_j(y) \leq 0$ for all $j \notin I(x)$, provided that $s \leq \delta_\eta$.

By virtue of implication (A.10), we are left with the task of proving the inequality $\theta(y) \leq \theta(x) + \lambda s \nabla \theta(x) F(x)$ for all $s \leq \delta_\eta$ and $x \in \widetilde{S}$ with $d(x) \geq \eta$. Using (A.6), we obtain for all $s \in [0, r]$:

$$\theta(y) \leq \theta(x) + s\nabla \theta(x) F(x) + \nabla \theta(x)(y - x - sF(x)) + K|y - x|^2 \quad (A.13)$$

where $K = \frac{1}{2}\max\left\{|\nabla^2 \theta(z)| : |z| \leq \beta\right\}$. The derivation of (A.13) follows from majorizing the second derivative of the mapping $w \to p(w) = \theta(x + w(y - x))$ and using the inequality $|y - x| \leq 2s|F(x)|$ (which is a direct consequence of (A.16)). It follows from (A.6), (A.5), (A.13) and (A.8) that the following inequality holds for all $s \in [0, r]$:

$$\theta(y) \leq \theta(x) + s\nabla \theta(x) F(x) + Qs^2\left(|\nabla \theta(x)||F(x)| + |F(x)|^2\right) \quad (A.14)$$

Definitions (A.3), (A.9) and inequality (A.14) allow us to conclude that $\theta(y) \leq \theta(x) + \lambda s \nabla \theta(x) F(x)$ for all $s \leq \delta_\eta$.

The proof is complete. ◁